\documentclass{amsart}
\usepackage{amsmath}
\usepackage{amsthm}
\usepackage{amssymb}
\usepackage{mathrsfs}
\usepackage{tikz}
\usetikzlibrary{matrix,arrows,positioning,calc}
\usetikzlibrary{shapes.multipart}
\usepackage{tikz-cd}
\pgfdeclarelayer{background layer}
\pgfsetlayers{background layer,main}

\newcommand{\cI}{{\mathcal I}}

 \newcommand{\cH}{{\mathscr H}}
 \newcommand{\cD}{{\mathscr D}}
 \newcommand{\cJ}{{\mathscr J}}
 \newcommand{\cR}{{\mathscr R}}
 \newcommand{\cL}{{\mathscr L}}

\newcommand{\Semigroups}{\textsc{Semigroups}}
\newcommand{\SgpDec}{\textsc{SgpDec}}
\newcommand{\GAP}{\textsc{Gap}}

\newcommand{\id}{1}
\newcommand{\quotient}{/\!\!}

\usepackage{url}
\usepackage{color}
\usepackage{colortbl}
\usepackage{hyperref}

\definecolor{darkblue}{rgb}{0.0,0.0,0.3}
\definecolor{grey}{rgb}{0.9,0.9,0.9}
\hypersetup{colorlinks,breaklinks,
            linkcolor=darkblue,urlcolor=darkblue,
            anchorcolor=darkblue,citecolor=darkblue}

\DeclareMathOperator{\im}{\lambda}

\newtheorem{theorem}{Theorem}[section]
\newtheorem{lemma}[theorem]{Lemma}
\newtheorem{corollary}[theorem]{Corollary}
\theoremstyle{definition}
\newtheorem{definition}[theorem]{Definition}
\newtheorem{example}{Example}
\newtheorem{fact}{Fact}

\begin{document}
\title[Skeleton Key]{Skeleton Key:  Subduction Classes in Finite Transformation Semigroups and Green's Relations}
\author{Attila Egri-Nagy
and Chrystopher L. Nehaniv}

\maketitle
\begin{abstract}
We establish key connections between Green's $\cJ$- and $\cL$-relations on a finite semigroup and the subduction relation defined on the image sets of an action of the same semigroup when it acts faithfully on a finite set.
The construction of the skeleton order, the partial order on equivalence classes of the subduction relation, is shown to depend in a functorial way on transformation semigroups and surjective morphisms, and to factor through the Green's $\leq_\cL$-order and  $\leq_\cJ$-order on the semigroup and  through the inclusion order on image sets.
For right regular representations, the correspondence between the $\cJ$-class order and the skeleton order is one of isomorphism. Finally, we characterize the relationship between natural subsystems of a transformation semigroup, permutator groups  and the $\cH$-relation.
\end{abstract}

\section{Introduction \& Mathematical Preliminaries}
The analysis of how a transformation semigroup acts on the subsets of the state set which occur as images under the semigroup action is a natural and fruitful direction of exploration.\footnote{It leads, for example, to   
the holonomy decomposition algorithm \cite{zeiger67a,zeiger68,ginzburg_book68,eilenberg,holcombe_textbook,KRTforCategories,automatanetworks2005} in Krohn-Rhodes theory, a widely applied wreath product decomposition theorem for finite transformation semigroups.} For a transformation semigroup, the study of this action on image sets is shown here depend on the fundamental Green's relations on the semigroup, and thus is of broad intrinsic interest in the theory of finite semigroups and their faithful actions. (See below for precise definitions.) 

One of the main tools of this analysis is the \emph{subduction} preorder relation defined on the set of images of the members of the semigroup considered as mappings.
 Green's preorders give ample information about the semigroup's internal structure, while subduction captures details of the semigroup action.
Therefore, the natural question arises:
\emph{What is the connection between the Green's relations and the subduction relation?}
More generally, by aiming to describe the connection between a semigroup action and the internal structure of semigroup itself we may yield insight and understanding into what transformation representations are possible for an abstract semigroup and the fine  details of their structure.  Lemma~\ref{lem:image-map} and the theorems proved in this paper are new, elucidating the key relationships of concepts well-known in classical semigroup theory to those used in holonomy methods, such as the skeleton. 

For clarity, we include some important mathematical preliminaries that will be
used throughout this manuscript:  
$(X,S)$ is a {\em transformation semigroup} with $S$ acting faithfully on the \emph{state set} $X$ if  $S$ is a subsemigroup of the (right) full transformation semigroup  ${\mathcal T}(X)$ of all mappings on a set $X$.
For $x\in X, s\in S$, the result of the action is written $x^s$. 
The action can be extended to subsets of $X$, if $P\subseteq X$ and $s\in S$ then $P^s=\{x^s\mid x\in P\}$.
The \emph{image} of a transformation $s$ is defined by $\im(s)=X^s$, and  we can also say that $\im(s)$ is the image of $X$ under $s$.
$S^\id$ is the monoid obtained by adjoining the identity on $X$ to the semigroup $S$, if it is not a member of $S$, otherwise $S^\id=S$.

$(A,\leq)$ is a \emph{preorder} (sometimes called a `quasi-order') if $\leq$ is a reflexive and transitive relation on the set $A$.
For a preorder, there exists an equivalence relation $(A,\equiv)$ defined by $a\equiv b \iff a\leq b$ and $b\leq a$,
and an induced partial order on the equivalence classes $(A/\!\equiv,\leq)$.
The surjective map $A\twoheadrightarrow A/\!\equiv$ is denoted by $\eta$.

The classical Green's relations $\leq_\cR$, $\leq_\cL$, $\leq_\cJ$ and  $\leq_\cH$, on any semigroup $S$ are the preorders:
$t \leq_\cR s \iff  tS^\id \subseteq sS^\id,$
$t \leq_\cL s \iff  S^\id t \subseteq S^\id s,$
$t \leq_\cJ s \iff  S^\id t S^\id \subseteq S^\id s S^\id$, and $\leq_\cH$ is the intersection of the $\leq_\cL$ and $\leq_\cR$ relations.
Then $ \cR,\cL,\cJ$ and $\cH$ denote the equivalence relations arising from the preorders  $\leq_\cR$, $\leq_\cL$, $\leq_\cJ$ and  $\leq_\cH$, respectively, and in each case the equivalence classes carry the induced partial order.   The equivalence relation $\cD$ on $S$ is the composite of $\cL$ and $\cR$, which commute. In the finite case, 
the $\cJ$ and $\cD$ relations coincide.
(For elementary properties of Green's preorders and their applications, see standard references for instance \cite{clifford_preston,Howie95,lallement1979}).
Here we only consider finite transformation semigroups.\footnote{
Nevertheless, we shall refer the $\cJ$-ordering and $\cJ$-classes, rather than  denoting $\cJ$- by $\cD$- as some authors working with finite semigroups do,  since  by definition $\cD$ is an equivalence relation that does {\em not} necessarily come from a preorder in the general setting of all (finite and infinite) semigroups  (see e.g.~\cite{lallement1979} for examples).}

\section{Subduction Relation}
For a transformation semigroup $(X,S)$ the set $\cI(X)=\{\im(s)\mid s\in S^\id\}$ is the \emph{image set} of the semigroup action. Note  $X=\im(\id)$ is always in $ \cI(X)$.

\begin{definition}[Subduction Pre-order]
Given a transformation semigroup $(X,S)$, for image sets $P$ and $Q$ in $\cI(X)$ 
$$P\subseteq_S  Q \iff \exists s \in S^\id  \text{ such that }  P \subseteq Q^s.$$
That is,  $P$ is a subset of $Q$, or it is a subset of some image of $Q$ under the semigroup action.
We say that $P$ is \emph{subduction below} of $Q$, or alternatively, $Q$ \emph{subduces} $P$.
\label{def:subduction}
\end{definition}
\begin{lemma}
 1. $(\cI(X),\subseteq_S)$ is a preorder. \\
 2. If $P \subseteq_S Q$ and $Q \subseteq_S P$ then $|P|=|Q|$.
\end{lemma}
\proof
Obviously $\subseteq_S$ is reflexive, since $P\subseteq P^\id$. It is transitive, since if  $P \subseteq Q^{s_1}$ and $ Q \subseteq R^{s_2}$ then $P \subseteq Q^{s_1} \subseteq (R^{s_2})^{s_1}=R^{s_2s_1}$, whence $P \subseteq_S R$.  To see (2),  there exists $s\in S^\id$ with $P \subseteq Q^s$,  so $Q$ has cardinality at least that of $P$. By symmetry, it follows that $P$ and $Q$ have the same cardinality. \qed\\

Therefore, we can naturally define by mutually subduced subsets an equivalence relation on $\cI(X)$, denoted by $\equiv_S$, leading to a partial order.\footnote{This structure provides the scaffolding for a holonomy decomposition since subduction equivalent subsets have isomorphic holonomy permutation groups, so only one copy of these groups is needed per class in the decomposition \cite{zeiger67a,zeiger68,ginzburg_book68,eilenberg,holcombe_textbook,automatanetworks2005}.
For the holonomy decomposition the skeleton order is extended by using
the {\em extended image set} $\cI^+(X)=\cI(X)\cup \{\{x\}: x \in X\}$  which includes any singletons that do not occur as images. This could potentially result in additional minimal equivalence classes for these singletons being adjoined to the skeleton.}

\begin{definition}[Skeleton Order]
The  {\em skeleton ordering}  $(\cI(X)\quotient\equiv_S,\subseteq_S)$ for transformation semigroup $(X,S)$ is the partial order on subduction equivalence classes of $\cI(X)$. 
\end{definition}

\section{$\cJ$-classes and Skeleton Classes}
We can establish connection between the induced classes of two preorders through a given preorder preserving map.
First, we describe the situation for preorders in general.

Let $(A_1,\leq_1)$ and $(A_2,\leq_2)$ be preorders.
Then, a function
$f:A_1\rightarrow A_2$ respects preordering, if
for all $a,a'\in A_1$, $a\leq_1 a'$ implies $f(a)\leq_2 f(a')$.

For a pre-order $(A, \leq)$, the equivalence class of $a \in A$, is $[a]=\{a' \in A : a\leq a' \mbox{ and } a'\leq a\}$.  The set of equivalence classes $A\quotient \equiv$  $\{[a]: a\in A\}$ carries the {\em induced partial order} $(A\quotient \equiv, \leq)$ given by $[a] \leq [a']$ if $a \leq a'$. 

The  {\em natural quotient map}  is the function  $\eta(a)=[a]$, which respects preordering and maps onto the induced partial order:
$\eta: (A, \leq) \rightarrow (A\quotient\equiv, \leq)$.
 
\begin{lemma}
\label{lem:preorders}
Given a function $f: (A_1,\leq_1) \rightarrow (A_2,\leq_2)$ that respects pre-ordering, we have that 
\begin{enumerate}
\item 
 $f$ induces an order-preserving map  $\bar{f}:(A_1\quotient\equiv_1,\leq_1)\rightarrow (A_2\quotient\equiv_2,\leq_2)$, and the following  diagram commutes,
 
\begin{center}
\begin{tikzcd}
A_1  \arrow{r}{f}\arrow[two heads]{d}{\eta_1} & A_2\arrow[two heads]{d}{\eta_2} \\
A_1\quotient\equiv_1 \arrow{r}{\bar{f}}& A_2\quotient\equiv_2
\end{tikzcd}
where $\eta_1$ and $\eta_2$ are the natural quotients
\end{center}
\item The equivalence relation induced on $A_1$ by the composite map
$\bar{f} \circ \eta_1  = \eta_2 \circ f$ is not finer than $\equiv_1$.
\end{enumerate}
\end{lemma}
\proof
(1) Let $\bar{f}$ denote the function taking the $\equiv_1$-class $[a]_1$ of any $a\in A_1$ to the $\equiv_2$-class $[f(a)]_2$ of $f(a)$.
If $a\equiv_1 a'$ then by definition $a\leq_1 a'$ and $a'\leq_1 a$.
Since $f$ respects preordering, $f(a)\leq_2 f(a')$ and $f(a')\leq_2 f(a)$, therefore $f(a)\equiv _2 f(a')$.
It follows that $\bar{f}$, given by $\bar{f}([a]_1)$=$[f(a)]_2$, is  well-defined and
order preserving.
(2) Since $\bar{f}$  is a function,  the inverse image $B$ of an $\equiv_2$-class is a set of  $\equiv_1$-classes. Hence the inverse image of $B$ in $A_1$ is
the union of these $\equiv_1$-equivalence classes. 
 \qed\\

\noindent{\bf Remark: } It is important to notice that $a <_1 a'$ does not exclude the possibility of $f(a)\equiv_2 f(a')$. Moreover,
even if  neither $a\leq_1 a'$ nor $a' \leq_1 a$ holds one might still have $f(a)<_2 f(a')$ or $f(a)\equiv_2 f(a')$. \\

Now we have two preorders: $\leq_\cJ$ on $S$ and subduction $\subseteq_S$ on $\cI(X)$.
Next we show that the surjective function $\im$ respects preordering.
For the weaker case, it is a basic fact that $a\cL b\implies \im(a)=\im(b)$.
However, $\cJ$-related elements can have different images.  For instance, in the full transformation semigroups on $n$ points, all constant maps are $\cJ$-equivalent to each other. %

\begin{lemma}\label{lem:image-map} For any transformation semigroup $(X,S)$ and
 any $a,b\in S$, we have
$$a\leq_\cL b\implies \im(a) \subseteq \im(b).$$
$$a\leq_\cJ b\implies \im(a)\subseteq_S\im(b).$$
That is,  $\im$ maps the $\cL$-preorder to the inclusion partial order and maps the $\cJ$-preorder to the subduction preorder.  Moreover, $\im$ induces a surjective map from $S^\id$ in each case.
\end{lemma}
\proof
The first assertion is well-known:
If $a \leq_\cL b$ then $a=sb$ for some $s\in S^\id$.  Thus $\im(a)=X^a=X^{sb}=(X^s)^b\subseteq X^b=\im(b)$.

For the second, if $a\leq_\cJ b$ then there exist $s,t\in S^\id$ such that $a=sbt$,
$$ \im(a)=\im(sbt)=\im(sb)^t\subseteq
\im(\id b)^t
=\im(b)^t,$$
therefore $\im(a)\subseteq_S\im(b)$.
Obviously $\im$ maps $S^1$ surjectively  onto $\cI(S)=\{\im(s): s \in S\}$, hence onto the preorder $(\cI(X),\subseteq_S)$ which has the same underlying set.
\qed

\begin{theorem}\label{functor}
 For a transformation semigroup $(X,S)$, there is a surjective order-preserving map $\bar{\im}_S$
from the partial order of $\cJ$-classes $({S^\id}/\cJ,\,\leq_{\cJ\!\!\!}~),$
onto the partial order of subduction classes $(\cI(X)\quotient\equiv_S,\subseteq_S)$.
The inverse image of a subduction equivalence class is a union of $\cJ$-classes.
\end{theorem}
\proof
$\im$ is a surjective, and is a preorder morphism from the Green's $\cJ$ preorder to the subduction preorder by Lemma~\ref{lem:image-map}, therefore by using Lemma~\ref{lem:preorders}(1), the induced map 
$\bar{\im}_S$ 
is a surjective order-preserving map. By Lemma~\ref{lem:preorders}(2), the inverse image of a subduction class corresponds to a union of $\cJ$-classes.~\qed\\

Similarly, generalizing the basic fact mentioned above, we have
\begin{theorem}
For a transformation semigroup $(X,S)$, there is a surjective order-preserving
map $\bar{\im}$ from the $\cL$-class order for $S^\id$ onto the inclusion partial order on $\cI(X)$.
The inverse image of an image set is a union of $\cL$-classes.
\end{theorem}

Putting these facts together, it follows that
\begin{theorem}\label{diagram}  For any transformation semigroup $(X,S)$, there  is a commutative diagram of surjective order-preserving morphisms:
\begin{center}
\begin{tikzcd}
(S^\id, \leq_\cL)\arrow[two heads]{d}{/\cL}\arrow[two heads]{rd}{\im}\\
(S^\id / \cL, \leq_\cL)   \arrow[two heads]{r}{\bar{\im}}\arrow[two heads]{d}{/\cJ} & (\cI(X),\subseteq) \arrow[two heads]{d}{\quotient\,\equiv_S} \\
 (S^\id / \cJ, \leq_\cJ)\arrow[two heads]{r}{\bar{\im}_S} & (\cI(X)\quotient \equiv_S,\subseteq_S)
\end{tikzcd}
\end{center}
\end{theorem}

\begin{corollary}
For the right regular representation $(S^\id,S)$:
\begin{enumerate}
\item The  $\cJ$-class order and the subduction order are isomorphic.
\item The $\cL$-class order and the inclusion order on image sets $\cI(X)$ are isomorphic.
\end{enumerate}
\end{corollary}
\proof
(1) By Lemma~\ref{lem:image-map}, it suffices to show that $\im(a)\subseteq_S\im(b)\implies a\leq_\cJ b$.
By definition of subduction
$\im(a)\subseteq\im(b)^t$  for some $t\in S^\id$.
Since $X=S^\id$ we can write $\im(a)$ as $\big(S^\id\big)^a$, or by shifting notation from semigroup action to semigroup multiplication, simply as $S^\id a$.
Therefore,
$$S^\id a\subseteq S^\id bt\implies S^\id a S^\id\subseteq S^\id b t S^\id\subseteq S^\id b S^\id\implies a\leq_\cJ b.$$  It follows that, if  $\im(a) \not\equiv_S \im(b)$ then $a \cJ b$ does not hold. Thus $\bar{\im}_S$ is injective, hence bijective.

\noindent (2)   More simply for the $\cL$-order, $\im(a) \subseteq \im(b)$ in the case of the right regular representation is just $(S^\id)^a \subseteq (S^\id)^b$, i.e., $S^\id a \subseteq S^\id b$,  which is the definition of
$a \leq_\cL b$.  Hence, $\im(a) \subseteq \im(b)$ implies $a\leq_\cL b$. By Lemma~\ref{lem:image-map} for the $\leq_{\cL}$-preorder, the converse holds. It follows that if $\im(a)\neq\im(b)$ then it cannot be that $a \cL b$, hence $\bar{\im}$ is injective, and hence bijective as well.
\qed\\

In the case of the right regular representation this says that the horizontal mappings in Theorem~\ref{diagram} are order isomorphisms.

Both the $\cJ$-class order and the skeleton capture information about the structure of the semigroup, therefore surjective homomorphisms should respect them.

\begin{theorem}[Functoriality]
Suppose $\varphi: (X,S) \twoheadrightarrow (Y,T)$ is a surjective morphism of transformation semigroups such that if $\id \in S$ then $\varphi(\id)$ is the identity on $Y$. Then $\varphi$ induces a natural mapping of the commutative
diagram for $(X,S)$ as in Theorem~\ref{diagram}, to the commutative diagram
for $(Y,T)$.
\end{theorem}
\proof A surjective map of semigroups induces a surjective map of the $\leq_\cL$ and $\leq_\cJ$ pre-orders and orderings (as well as for $\leq_\cR$ and $\leq_\cH$).
$\varphi$ also induces a surjective map from $\cI(X)$ onto $\cI(Y)$, and subduction in the source implies subduction in the target since $P \subseteq Q^s$ implies
$\varphi(P) \subseteq \varphi(Q^s)=\varphi(Q)^{\varphi(s)}$, hence the subduction relation is respected, and the result follows. 
  \qed

\vspace{-.5cm}
\section{Examples}

We present a few examples to illustrate the connection between the $\cJ$-class order and the skeleton. The partial orders are displayed as Hasse diagrams. Shaded clusters of $\cJ$-classes are mapped to a single subduction class.
\begin{example}[Simple collapsing of a chain]  Let $X=\{1,2,3\}$,
$t_1=\left(\begin{smallmatrix}1&2&3\\1& 3& 3\end{smallmatrix}\right)$,
$t_2=\left(\begin{smallmatrix}1&2&3\\3& 1& 3\end{smallmatrix}\right)$,
$t_3=\left(\begin{smallmatrix}1&2&3\\3& 3& 3\end{smallmatrix}\right)$
and $M$ the monoid $\{\id,t_1,t_2,t_3\}$, so $(X,M)$ is a transformation monoid on 3 points. The principal two-sided ideals are:
$M\id M=M$, $Mt_1M=\{t_1,t_2,t_3\}$, $Mt_2M=\{t_2,t_3\}$, $Mt_3M=\{t_3\}$,
therefore $t_3 <_{\cJ} t_2 <_{\cJ} t_1 <_{\cJ} \id$ and all elements form a singleton $\cJ$-class on their own. $\cI(X)=\big\{\{1,2,3\},\{1,3\},\{3\}\big\}$ defines the subduction classes.
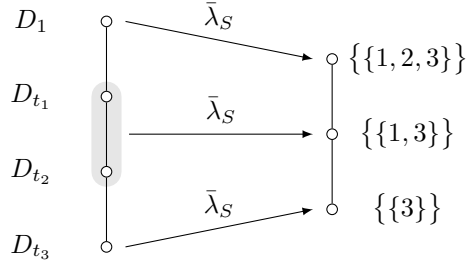
\begin{center}
\begin{tikzpicture}
  \tikzstyle{every node}=[draw,circle,fill=white,minimum size=4pt,inner sep=0pt]
  \tikzstyle{label}=[draw=none,fill=none,inner sep=1pt]
\draw[line width=0.4cm,color=black!10,cap=round,join=round] (0,1)--(0,2);

\draw (0,0) node (Dt3) {};
\draw (0,1) node (Dt2) {};
\draw (0,2) node (Dt1) {};
\draw (0,3) node (D1) {};

\draw node [left of=Dt3,label] {$D_{t_3}$} ;
\draw node [left of=Dt2,label] {$D_{t_2}$} ;
\draw node [left of=Dt1,label] {$D_{t_1}$} ;
\draw node [left of=D1,label] {$D_{1}$} ;

\draw (Dt3) -- (Dt2);
\draw (Dt2) -- (Dt1);
\draw (Dt1) -- (D1);

\draw (3,0.5) node (3) {};
\draw (3,1.5) node (13) {};
\draw (3,2.5) node (123) {};

\draw (3) -- (13);
\draw (13) -- (123);

\draw node [right of=3,label] {$\big\{\{3\}\big\}$} ;
\draw node [right of=13,label] {$\big\{\{1,3\}\big\}$} ;
\draw node [right of=123,label] {$\big\{\{1,2,3\}\big\}$} ;

 \draw[->,>=latex] ($(D1.east)+(5pt,0pt)$) to node[above,midway,label]{$\bar{\im}_S$} ($(123.west)+(-5pt,0pt)$) ; 
 \draw[->,>=latex] ($(Dt3.east)+(5pt,0pt)$) to node[above,midway,label]{$\bar{\im}_S$} ($(3.west)+(-5pt,0pt)$) ; 
 \draw[->,>=latex] (.3,1.5) to node[above,midway,label]{$\bar{\im}_S$} ($(13.west)+(-5pt,0pt)$) ; 

\end{tikzpicture}
\end{center}
\end{example}
A simple linear order is mapped to a shorter linear order, since $\im(t_1)=\im(t_2)=\{1,3\}$.

\begin{example} More general collapsing (a usual motif) for a transformation monoid on 3 points,
$M=\big\{\id,
\left(\begin{smallmatrix}1&2&3\\1& 1& 3\end{smallmatrix}\right),
\left(\begin{smallmatrix}1&2&3\\3& 2& 3\end{smallmatrix}\right),
\left(\begin{smallmatrix}1&2&3\\3& 1& 3\end{smallmatrix}\right),
\left(\begin{smallmatrix}1&2&3\\3& 3& 3\end{smallmatrix}\right)\big\}.$
The right regular transformation representation of $M$ can be encoded as
$$M'=\big\{\id,
\left(\begin{smallmatrix}1&2&3&4&5\\2& 2& 4& 4& 5\end{smallmatrix}\right),
\left(\begin{smallmatrix}1&2&3&4&5\\3& 5& 3& 5& 5\end{smallmatrix}\right),
\left(\begin{smallmatrix}1&2&3&4&5\\4& 5& 4& 5& 5\end{smallmatrix}\right),
\left(\begin{smallmatrix}1&2&3&4&5\\5& 5& 5& 5& 5\end{smallmatrix}\right)\big\}.$$  
Its skeleton is isomorphic to the $\cJ$-class order of $M$.
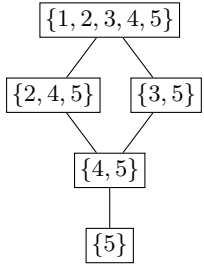
\begin{center}
  \begin{tikzpicture}
  \tikzstyle{every node}=[draw,shape=rectangle,inner sep=2pt]
  \tikzstyle{label}=[draw,fill=none,inner sep=1pt]
\small

\draw (0,0) node (D1) {$\{5\}$};
\draw (0,1) node (D2) {$\{4,5\}$};
\draw (-.75,2) node (D3) {$\{2,4,5\}$};
\draw (.75,2) node (D4) {$\{3,5\}$};
\draw (0,3) node (D5) {$\{1,2,3,4,5\}$};

\draw (D1) -- (D2);
\draw (D2) -- (D3);
\draw (D2) -- (D4);
\draw (D3) -- (D5);
\draw (D4) -- (D5);

\end{tikzpicture}
  \begin{tikzpicture}
  \tikzstyle{every node}=[draw,circle,fill=white,minimum size=4pt,inner sep=0pt]
  \tikzstyle{label}=[draw=none,fill=none,inner sep=1pt]
\draw[line width=0.4cm,color=black!10,cap=round,join=round] (0,1)--(.75,2);

\draw (0,0) node (D1) {};
\draw (0,1) node (D2) {};
\draw (-.75,2) node (D3) {};
\draw (.75,2) node (D4) {};
\draw (0,3) node (D5) {};

\draw (D1) -- (D2);
\draw (D2) -- (D3);
\draw (D2) -- (D4);
\draw (D3) -- (D5);
\draw (D4) -- (D5);

\draw (3,0) node (S1) {};
\draw (3,1) node (S2) {};
\draw (3,2) node (S3) {};
\draw (3,3) node (S4) {};

\draw (S1) -- (S2);
\draw (S2) -- (S3);
\draw (S3) -- (S4);

 \draw[->,>=latex] (.6,1.5) to node[above,midway,label]{$\bar{\im}_S$} ($(S2.west)+(-5pt,0pt)$) ; 
\end{tikzpicture}
\end{center}
\end{example}


\begin{example}
$M$ monoid  generated by
$a=\left(\begin{smallmatrix}1&2&3&4&5\\2& 1& 1& 1&4\end{smallmatrix}\right)$ and
$b=\left(\begin{smallmatrix}1&2&3&4&5\\1& 2& 2& 3&4\end{smallmatrix}\right)$.
In $M$,  $a$ and $b$ are $\leq_\cJ$-incomparable, but $\im(a)\subset_S\im(b)$. This is so as  there is no solution for the equation $b=sat$ or $a=sbt$ for $s,t\in M$, although $\im(a)\subset\im(b)$.

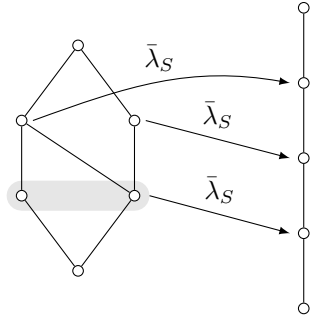
\begin{center}
\begin{tikzpicture}
  \tikzstyle{every node}=[draw,circle,fill=white,minimum size=4pt,inner sep=0pt]
  \tikzstyle{label}=[draw=none,fill=none,inner sep=1pt]

\draw[line width=0.4cm,color=black!10,cap=round,join=round] (-.75,1)--(.75,1);

\draw (0,0) node (D1) {};
\draw (-.75,1) node (D2) {};
\draw (.75,1) node (D3) {};
\draw (-.75,2) node (D4) {};
\draw (.75,2) node (D5) {};
\draw (0,3) node (D6) {};

\draw (D1) -- (D2);
\draw (D1) -- (D3);
\draw (D2) -- (D4);
\draw (D3) -- (D4);
\draw (D3) -- (D5);
\draw (D4) -- (D6);
\draw (D5) -- (D6);

\draw (3,-0.5) node (S1) {};
\draw (3,.5) node (S2) {};
\draw (3,1.5) node (S3) {};
\draw (3,2.5) node (S4) {};
\draw (3,3.5) node (S5) {};

\draw (S1) -- (S2);
\draw (S2) -- (S3);
\draw (S3) -- (S4);
\draw (S4) -- (S5);

 \draw[->,>=latex] (.95,1) to node[above,midway,label]{$\bar{\im}_S$} ($(S2.west)+(-3pt,0pt)$) ; 
 \draw[->,>=latex] (.9,2) to node[above,midway,label]{$\bar{\im}_S$} ($(S3.west)+(-3pt,0pt)$) ; 
 \draw[->,>=latex,out=20,in=170] (-.6,2) to node[above,midway,label]{$\bar{\im}_S$} ($(S4.west)+(-3pt,0pt)$) ; 

\end{tikzpicture}
\end{center}
This shows that the subduction order may contain new relations beyond those induced by collapsing nodes of the $\cJ$-order diagram.
Consequently, the length of a longest $\cJ$-chain is not an upper bound for the height of the skeleton.
\end{example}

So far the $\cJ$-class orders were all lattices, but this is not true in general, therefore we have to look at a monoid with more inner structure.

\begin{example}[Nonlinear, non-lattice skeleton] Let
$a=\left(\begin{smallmatrix}1&2&3&4&5\\2& 2& 1& 2&4\end{smallmatrix}\right)$,
$b=\left(\begin{smallmatrix}1&2&3&4&5\\3& 5& 2& 3&2\end{smallmatrix}\right)$,
$c=\left(\begin{smallmatrix}1&2&3&4&5\\3& 5& 4& 5&4\end{smallmatrix}\right)$ and $M=\langle a,b,c\rangle$.
$|M|=31$, $|\cI(X)|=16$, number of $\cD$-classes is 13, and the number of skeleton classes is 9.
On the left the $\cD$-class picture is drawn. On top of each $\cL$-class (drawn vertically) the corresponding image is displayed.
$\cH$-classes with an idempotent are shaded.
The grey background blobs indicate $\cD$-classes that are collapsed into one subduction class.
On the right the skeleton order is drawn.
It is nonlinear and it is not a lattice. The boxes indicate subduction equivalence classes.

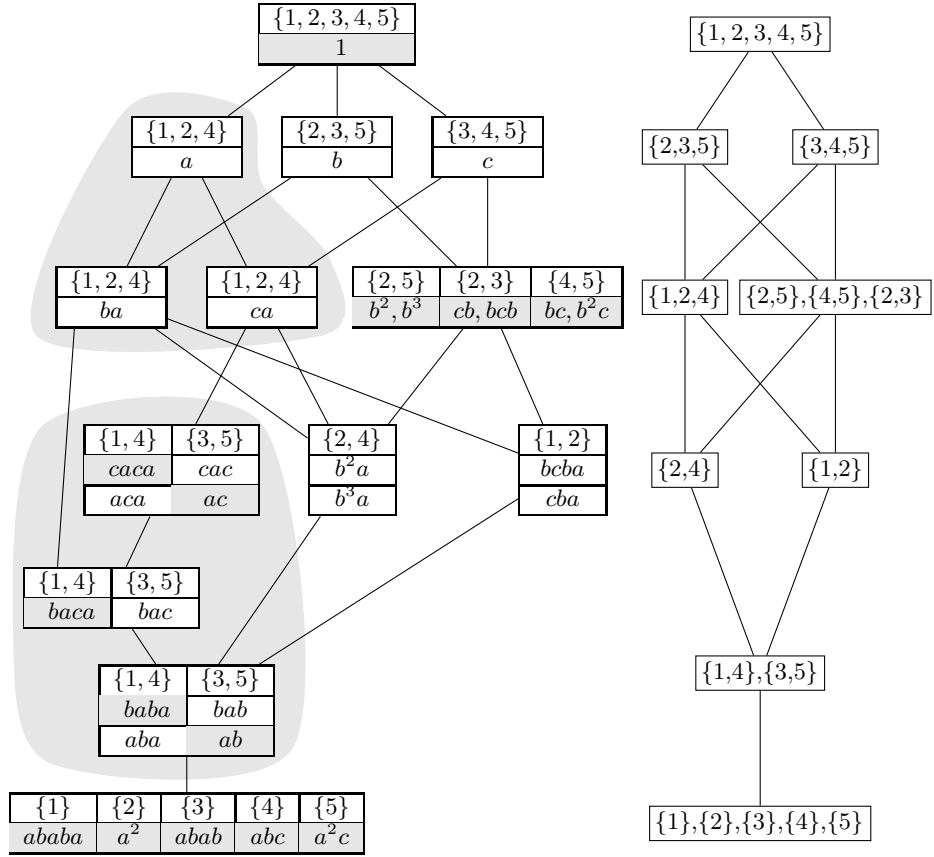
\begin{center}
\begin{tikzpicture}
[
  tab/.style={rectangle,inner sep=0pt},
  mySingleItem/.style={rectangle,draw=black,thin}]  
\small
 \node[tab] at (0,0) (n1) 
    {\begin{tabular}{|c|}\hline
       $\{1,2,3,4,5\}$ \\ \hline
        \cellcolor{grey} 1\\ \hline 
    \end{tabular}
    };

 \node[tab,fill=white] at (-2,-1.5) (na)
    {\begin{tabular}{|c|}\hline
       $\{1,2,4\}$ \\ \hline
        $a$\\ \hline 
    \end{tabular}
    };

 \node[tab] at (0,-1.5) (nb)
    {\begin{tabular}{|c|}\hline
       $\{2,3,5\}$ \\ \hline
        $b$\\ \hline 
    \end{tabular}
    };

 \node[tab]  at (2,-1.5) (nc)
    {\begin{tabular}{|c|}\hline
       $\{3,4,5\}$ \\ \hline
        $c$\\ \hline
    \end{tabular}
    };

 \node[tab,fill=white]  at (-3,-3.5) (nba)
    {\begin{tabular}{|c|}\hline
       $\{1,2,4\}$ \\ \hline
        $ba$\\ \hline 
    \end{tabular}
    };

 \node[tab,fill=white]  at (-1,-3.5) (nca)
    {\begin{tabular}{|c|}\hline
       $\{1,2,4\}$ \\ \hline
        $ca$\\ \hline 
    \end{tabular}
    };

 \node[tab] at (2,-3.5) (nb2)
    {\begin{tabular}{|c|c|c|}\hline
       $\{2,5\}$&$\{2,3\}$&$\{4,5\}$ \\ \hline
        \rowcolor{grey} $b^2,b^3$ & $cb,bcb$ & $bc,b^2c$\\ \hline 
    \end{tabular}
    };

 \node[tab,fill=white] at (-2.2,-5.8) (ncaca)
    {\begin{tabular}{|c|c|}\hline
       $\{1,4\}$&$\{3,5\}$ \\ \hline
        \cellcolor{grey} $caca$ & $cac$ \\ \hline 
        $aca$ & \cellcolor{grey}$ac$ \\ \hline 
    \end{tabular}
    };

 \node[tab]  at (0.2,-5.8) (nb2a)
    {\begin{tabular}{|c|}\hline
       $\{2,4\}$ \\ \hline
        $b^2a$\\ \hline 
        $b^3a$\\ \hline 
    \end{tabular}
    };

 \node[tab]  at (3,-5.8) (nbcba)
    {\begin{tabular}{|c|}\hline
       $\{1,2\}$ \\ \hline
        $bcba$\\ \hline 
        $cba$\\ \hline 
    \end{tabular}
    };

 \node[tab,fill=white] at (-3,-7.5) (nbaca)
    {\begin{tabular}{|c|c|}\hline
       $\{1,4\}$&$\{3,5\}$ \\ \hline
        \cellcolor{grey} $baca$ & $bac$ \\ \hline 
    \end{tabular}
    };

 \node[tab,fill=white] at (-2,-9) (nbaba)
    {\begin{tabular}{|c|c|}\hline
       $\{1,4\}$&$\{3,5\}$ \\ \hline
        \cellcolor{grey} $baba$ & $bab$ \\ \hline 
        $aba$ & \cellcolor{grey}$ab$ \\ \hline 
    \end{tabular}
    };

 \node[tab] at (-2,-10.5) (nababa)
    {\begin{tabular}{|c|c|c|c|c|}\hline
       $\{1\}$&$\{2\}$&$\{3\}$&$\{4\}$&$\{5\}$ \\ \hline
        \rowcolor{grey} $ababa$ & $a^2$ & $abab$ & $abc$ & $a^2c$\\ \hline 
    \end{tabular}
    };

\draw (n1) -- (na);
\draw (n1) -- (nb);
\draw (n1) -- (nc);
\draw (na) -- (nba);
\draw (na) -- (nca);
\draw (nb) -- (nba);
\draw (nb) -- (nb2);
\draw (nc) -- (nca);
\draw (nc) -- (nb2);
\draw (nba) -- (nb2a);
\draw (nba) -- (nbcba);
\draw (nca) -- (ncaca);
\draw (nca) -- (nb2a);
\draw (nb2) -- (nbcba);
\draw (nb2) -- (nb2a);
\draw (nb2a) -- (nbaba);
\draw (nbcba) -- (nbaba);
\path  (nba.220) edge (nbaca.150);
\draw(ncaca) -- (nbaca);
\draw(nbaca) -- (nbaba);
\draw(nbaba) -- (nababa);
\begin{pgfonlayer}{background layer}
\filldraw [grey] plot [smooth cycle] coordinates {(-2.8,-1)(-1.2,-1) (-1,-2) (-0.1,-3)(-0.2,-4) (-4,-4)};
\filldraw [grey] plot [smooth cycle] coordinates {(-3.9,-5.2) (-1,-5.2) (-0.6,-9.5)(-3.8,-9.7)(-4.3,-8)};
\end{pgfonlayer}
\begin{scope}[xshift=160pt,node distance=50pt]
  \tikzstyle{every node}=[draw,shape=rectangle,inner sep=2pt]
  \tikzstyle{label}=[draw,fill=none,inner sep=1pt]
\small

\draw  node at(0,0) (n12345) {$\{1,2,3,4,5\}$};
\draw node at (-1,-1.5) (n235) {\{2,3,5\}};
\draw node at (1,-1.5) (n345) {\{3,4,5\}};
\draw node at (-1,-3.5)  (n124) {\{1,2,4\}};
\draw node at (1,-3.5) (n25) {\{2,5\},\{4,5\},\{2,3\}};
\draw node at (-1,-5.8) (n24) {\{2,4\}};
\draw node at (1,-5.8) (n12) {\{1,2\}};
\draw node at (0,-8.5)  (n14) {\{1,4\},\{3,5\}};
\draw node at (0,-10.5) (n1) {\{1\},\{2\},\{3\},\{4\},\{5\}};

\draw (n12345) -- (n235);
\draw (n12345) -- (n345);
\draw (n345) -- (n124);
\draw (n235) -- (n124);
\draw (n345) -- (n25);
\draw (n235) -- (n25);
\draw (n124) -- (n24);
\draw (n124) -- (n12);
\draw (n25) -- (n24);
\draw (n25) -- (n12);
\draw (n24) -- (n14);
\draw (n12) -- (n14);
\draw (n14) -- (n1);

\end{scope}
\end{tikzpicture}

\end{center}
The skeleton also contains nonsingleton subduction equivalence classes.
\end{example}

\section{Conclusion}
We clarified the connection between the $\cJ$-classes of a semigroup and the subduction classes of a transformation representation of the same semigroup.
We showed how the partial order of $\cJ$-classes constrains the image relations in the possible (faithful) actions of the semigroup.
Therefore, these results may also be useful for investigating or enumerating the possible action representations of a semigroup, and for investigating structures used in decompositions that are respected under surjective maps. 

For calculating and checking the examples we used the \GAP~\cite{GAP4} packages \Semigroups~\cite{Semigroups}, \SgpDec~\cite{SgpDec} and \textsc{SgpViz}\cite{sgpviz}.

\section*{Appendix: Natural Subsystems, Permutators, $\cH$-classes}
If $s\in S^1$ and $P \subseteq P^s$, then $P^s=P$ by finiteness and so $s$ permutes $P$.
The {\em permutator} of $P$ denoted $Perm(P)$ is the group of all such permutations of $P$, and consists of restrictions of all such $s$ to $P$. 
If $P=X^e$ for an idempotent $e$, and  $G_e$ is  the unique largest subgroup of $S$
containing $e$, then $(P,G_e)$ is called a {\em natural subsystem} of $(X,S)$.

\begin{theorem} \label{HNatSysThm} 
A natural subsystem $(P,G_e)$ is a faithful permutation group,~and
\begin{enumerate}
\item  \label{G_eHasAllPerms}
$(P,G_e)\cong (P, Perm(P))$, where the isomorphism from $G_e$ to $Perm(P)$
is given by restriction to $P$.  Moreover $G_e=H(e)$, the $\cH$-class of $e$.
\item \label{NatSystemsIso}
Also, for $\cJ$-equivalent idempotents $e_1$ and $e_2$,  $(X^{e_1},G_{e_1})\cong (X^{e_2},G_{e_2})$, 
and for $\cL$-equivalent idempotents, the state sets coincide. 
\end{enumerate}
\end{theorem}

The proof relies on some facts:

\begin{fact}\label{Fact1}
If $e^2=e$ and $s\in S$ permutes $X^e$, then $s$, $ese$, $es$, and $se$ agree as transformations on $X^e$.
\end{fact}
\proof
Take any $z=x\cdot e\in X^e$. Then $e$ acts as the identity on $X^e$, since $z\cdot e = (x\cdot e)\cdot e=x\cdot e^2=z$. Since $s$ maps $X^e$ to $X^e$, it follows that 
$z\cdot es= (z\cdot es)\cdot e= z\cdot ese = (z\cdot e) \cdot se= z\cdot se=(z\cdot s)\cdot e = z\cdot s.$
\qed

\begin{fact}\label{Fact2}
For $e^2=e$ and $s\in S$ permutes $P$.  Let $w>1$ such that $(ese)^w$ is the unique idempotent power of $ese$.  Then (1) $(ese)^w=e$ and  (2) $ese \in G_e$.  
\end{fact}
\proof
By Fact~\ref{Fact1}, $ese$ permutes $X^e$ and agrees with $s$ on $X^e$. 
Obviously $X^{ese}=(X^{es})^e \subseteq X^e$.   Now $(ese)^w$ is idempotent and so fixes its image, which is $X^e$, since $e$ fixes $X^e$, $s$ permutes $X^e$ and then $e$ fixes the result, for $w$ iterations.
Take $x\in X$. Now,
$x\cdot (ese)^w= (x\cdot e) \cdot  (ese)^w = x\cdot e$.  Thus the two transformations $e$ and $(ese)^w$ agree for all $x\in X$.  Since $(X,S)$ is faithful, $s =(ese)^w$.   
Since $w>1$, $(ese)^{w-1}(ese)=(ese)^w=(ese)(ese)^{w-1}$, and as the element in middle is $e$, we have shown that $(ese)^{w-1}$ is a multiplicative inverse of $ese$ with respect to $e$.  Therefore $ese$ is contained in a group $\langle ese \rangle$ having two-sided identity element $e$. Therefore $ese \in G_e$, the unique maximal subgroup of $S$ containing $e$. 
\qed\\

\noindent{\bf \it Proof of Theorem~\ref{HNatSysThm}.} 
That $(P,G_e)$ is a permutation group is is well-known and can be shown as follows: 
To see that $(P,G_e)$ is faithful, suppose that $g,g' \in G_e$ have the same action on all of $P=X^e$. Let $x \in X$, then $x\cdot g = x\cdot eg = 
(x\cdot e)\cdot g = (x\cdot e)\cdot g' = x\cdot g'$, since $x\cdot e \in P$. Thus $g$ and $g'$ have the same action on $X$, whence $g = g'$ since $(X, S)$ is faithful. Since $e$ acts as the identity on $P=X^e$, it is easy to check that $G_e$ acts on $P$ by permutations, with each $g^{-1}$ acting as the inverse of the mapping on $X^e$ due to $g\in G_e$.

Facts~\ref{Fact1} and \ref{Fact2} together show that any permutation of $X^e$ realized by some $s\in S$ can be realized by a member of $G_e$, namely $ese$.
That is, for all elements $s$ of $S$ permuting $X^e$, there is a transformation among the elements of the single group $G_e$ having the same restriction to domain $X^e$. Indeed, $s\restriction_{X^e}=ese\restriction_{X^e}$. The
maximal subgroup $G_e$ of $S$ with identity $e$  is $H(e)$, the $\cH$-class of an idempotent (cf.\ e.g.~\cite{Howie95}). This proves Theorem~\ref{HNatSysThm}(\ref{G_eHasAllPerms}), while Theorem~\ref{HNatSysThm}(\ref{NatSystemsIso}) follows from Lemma~\ref{lem:image-map} that $\cJ$-equivalent elements are subduction equivalent and $\cL$-equivalent elements have the same image. 
Since subduction equivalent image sets have isomorphic permutator groups (e.g. \cite{zeiger67a,eilenberg}),
the proof is complete.  \qed\\

\noindent
{\bf Acknowledgements.} This research was supported in part by the Natural Sciences and Engineering Research Council of Canada (NSERC), funding reference number RGPIN-2019-04669. 
Cette recherche a \'et\'e financ\'ee en partie par le Conseil de recherches en sciences naturelles et en g\'enie du Canada (CRSNG), num\'ero de r\'ef\'erence RGPIN-2019-04669.

\bibliographystyle{plain}
\bibliography{../coords}

\begin{thebibliography}{10}

\bibitem{clifford_preston}
A.H. Clifford and G.B. Preston.
\newblock {\em The Algebraic Theory of Semigroups, Vol.~1}.
\newblock Number~7 in Mathematical Surveys. American Mathematical Society, 2nd
  edition, 1967.

\bibitem{sgpviz}
Manuel Delgado and Jos{\'e} Morais.
\newblock {\em {{GAP} package \textsc{SgpViz}} 0.999.5}, 2022.
\newblock
  \href{https://gap-packages.github.io/sgpviz/}{\url{https://gap-packages.github.io/sgpviz/}}.

\bibitem{automatanetworks2005}
P{\'a}l D{\"o}m{\"o}si and Chrystopher~L. Nehaniv.
\newblock {\em {Algebraic Theory of Finite Automata Networks: An
  Introduction}}, volume~11 of {\em {SIAM Series on Discrete Mathematics and
  Applications}}.
\newblock Society for Industrial and Applied Mathematics, 2005.

\bibitem{SgpDec}
A.~Egri-Nagy, C.~L. Nehaniv, and J.~D. Mitchell.
\newblock {\em {\textsc{{S}gp{D}ec} -- software package for Hierarchical
  Composition and Decomposition of Permutation Groups and Transformation
  Semigroups, Version 1.1.0}}, 2024.
\newblock \url{https://gap-packages.github.io/sgpdec/}.

\bibitem{eilenberg}
Samuel Eilenberg.
\newblock {\em {Automata, Languages and Machines, vol.\ B}}.
\newblock Academic Press, 1976.

\bibitem{GAP4}
The GAP~Group.
\newblock {\em {GAP -- Groups, Algorithms, and Programming, Ver.\ 4.14.0}},
  2024.

\bibitem{ginzburg_book68}
Abraham Ginzburg.
\newblock {\em {Algebraic Theory of Automata}}.
\newblock Academic Press, 1968.

\bibitem{holcombe_textbook}
W.~M.~L. Holcombe.
\newblock {\em {Algebraic Automata Theory}}.
\newblock Cambridge University Press, 1982.

\bibitem{Howie95}
John~M. Howie.
\newblock {\em {Fundamentals of Semigroup Theory}}, volume~12 of {\em {London
  Mathematical Society Monographs New Series}}.
\newblock Oxford University Press, 1995.

\bibitem{lallement1979}
Gerard Lallement.
\newblock {\em Semigroups and Combinatorial Applications}.
\newblock Pure and Applied Mathematics. Wiley, New York, 1976.

\bibitem{Semigroups}
J.~D. Mitchell et~al.
\newblock {\em Semigroups - GAP package, Version 5.3.7}, Mar 2024.

\bibitem{KRTforCategories}
Charles Wells.
\newblock A {K}rohn-{R}hodes theorem for categories.
\newblock {\em Journal of Algebra}, 64:37--45, 1980.

\bibitem{zeiger67a}
H.~Paul Zeiger.
\newblock {Cascade synthesis of finite state machines}.
\newblock {\em Information and Control}, 10(4):419--433, 1967.
\newblock Erratum: {\bf 11}(4): 471 (1967).

\bibitem{zeiger68}
H.~Paul Zeiger.
\newblock {Cascade Decomposition Using Covers}.
\newblock In Michael~A. Arbib, editor, {\em {Algebraic Theory of Machines,
  Languages, and Semigroups}}, chapter~4, pages 55--80. Academic Press, 1968.

\end{thebibliography}
\end{document}